# On Asymptotic Expansion in the Random Allocation of Particles by Sets

Saidbek Sh. Mirakhmedov[1] and Sherzod M. Mirakhmnedov[2]

**Abstract.**

We consider a scheme of equiprobable allocation of particles into cells by sets. The Edgeworth type asymptotic expansion in the local central limit theorem for a number of empty cells left after allocation of all sets of particles is derived.

1. **Introduction.**

Many combinatorial problems in probability and statistics can be formulated and best understood by using appropriate random allocation schemes (alternatively known as urn models). Such models naturally arise in statistical mechanics (Feller (1968)), clinical trials (Gani (1993)), cryptography (Menezes et al. (1997)) etc. The properties of various random allocation schemes have been extensively studied in the probabilistic and statistical literature, see books by Kolchin et al. (1976), Johnson and Kotz (1977) and survey papers by Ivanov et al. (1984) and Kotz and Balakrishnan (1997).

The classical allocation scheme assumes equiprobable allocation of *n* particles to a finite number *N* of cells, i.e. the probability of a particle falling into any particular cell is equal to $N^{-1}$. There are several generalizations of the classical scheme, see, for example, Ivanov et al. (1984). In the present paper we consider the following generalization of the classical model.

Let cells be indexed by $1, 2, ..., N$ and sets of particles be indexed by $1, 2, ..., s$, a set of particles with index *j* containing $n_j$ particles. We assume that

- Particles are allocated one set at a time
- Sets of particles are allocated independently of each other.
- No two particles from a given set can be allocated to a common cell.

[1] Saidbek Sh. Mirakhmedov

Institute of "Algoritm and Engineering". Fayzulla Hodjaev-45, Tashkent. Uzbekistan

e-mail saeed_0810@yahoo.com

[2] Sherzod M. Mirakhmedov (former Sh.A. Mirakhmedov)(✉) GIK Institute of Engineering Sciences and Technology. Topi-23460. Swabi, N.W.F.P. Pakistan.

e-mail shmirakhmedov@yahoo.com ;

- All $C_N^{n_i}$ possible variants of the allocation of particles from the $i$-th set have the same probability equal to $\left(C_N^{n_i}\right)^{-1}$, $i=1,2,...,s$.

This paper deals with the random variable (r.v.) $\mu_0 = \mu_0(N, n_1,...,n_s)$ giving the number of empty cells left after all the $s$ sets of particles have been allocated.

Note that this allocation scheme corresponds to an urn model when we deal with $s$ independent samples of sizes $n_1, n_2,...,n_s$ respectively, where each sample is drawn by simple sample scheme without replacement from a population of size $N$; here the r.v. $\mu_0$ is a number of untouched elements of the population.

The classical allocation scheme corresponds to the case $n_1 = n_2 = ... = n_s = 1$. The random allocation of particles by sets was also studied in the literature. For detailed surveys of theory and applications see Ivanov et al. (1984) and Kotz and Balakrishman (1997). The most relevant to the present paper references are Park (1981), Mikhhaylov (1980, 1981), Vatutin and Mikhaylov (1982) and Mirakhmedov (1985, 1994, 1996). In particular, the most general limit theorem for r.v. $\mu_0$ were obtained by Mikhaylov and Vatutin (1982): they proved the central limit theorem and Poisson limit theorem for $\mu_0$ under rather weak conditions, giving lower bounds for the remainder terms. But the topics of interest to us here are not readily available in the literature. In the present paper, we will derive an Edgeworth type asymptotic expansion in the local central limit theorem for $\mu_0$ using a combination of two approaches. The first approach is based on a Bartlett type formula for the characteristic function of $\mu_0$ and the second one is based on the Frobenius-Harper technique, exploiting the fact that the r.v. $\mu_0$ can be represented as a sum of independent Bernoulli r.v.'s.

Thus in this paper our purpose is twofold. On the one hand, we wish to derive a better local approximation formula for the distribution of $\mu_0$. On the other hand, we would also like to develop the above mentioned methods to obtain an asymptotical expansion type result.

The organization of the paper is as follows. In Section 2, a Bartlett type formula is presented. An asymptotic expansion of the characteristic function of $\mu_0$ is given in Section 3. In Section 4, an asymptotic expansion is obtained in the central local limit theorem for $\mu_0$ in the case of a fixed $s$. The case when $s$ can increase as $\min_{1 \le l \le s} n_l \to \infty$ is discussed in Section 5.

In what follows, all the asymptotic relations assume that $\min_{1\leq l\leq s} n_l \to \infty$ and $N \to \infty$. By the symbol $C$ (with and without subscripts) we denote absolute positive constant. The relation $\xi \sim B(1, p)$ means that the r.v. $\xi$ has the Bernoulli distribution with parameter $p$, $0 < p < 1$.

## 2. Bartlett's type formula.

The allocation of the $l$-th set of particles is determined by the random allocation indicators $\eta_{l,m}$, $m = 1, 2, ..., N$, where $\eta_{l,m} = 1$ if a particle of the $l$-th set falls into the cell with index $m$, $\eta_{l,m} = 0$ otherwise. It is clear that $\eta_{l,m} \sim Bi(1, p_l)$ with $p_l = n_l / N$; also $\eta_{l,1} + \eta_{l,2} + ... + \eta_{l,N} = n_l$, $l = 1, 2, ..., s$, and

$$\mu_0 = \sum_{m=1}^{N} \mathbb{1}\{\eta_{1,m} + \eta_{2,m} + ... + \eta_{s,m} = 0\}, \tag{2.1}$$

where $\mathbb{1}\{A\}$ is the indicator of the event $A$.

Let $\xi_{l,m}$, $l = 1, 2, ..., s$; $m = 1, 2, ..., N$ be a collection of the mutually independent r.v.'s, $\xi_{lm} \sim Bi(1, p_l)$, $\zeta_{l,N} = \xi_{l,1} + \xi_{l,2} + ... + \xi_{l,N}$. It is easy to check that for any $k_1, k_2, ..., k_N$ with $k_l = 0$ or $1$ such that $k_1 + k_2 + ... + k_N = n_l$

$$P\{\eta_{l,1} = k_1, \eta_{l,2} = k_2 ..., \eta_{l,N} = k_N\} = P\{\xi_{l,1} = k_1, \xi_{l,2} = k_2 ..., \xi_{l,N} = k_N / \zeta_{l,N} = n_l\}. \tag{2.2}$$

Let $\bar{\eta}_l = (\eta_{l,1}, \eta_{l,2}, ..., \eta_{l,N})$, $\bar{\xi}_l = (\xi_{l,1}, \xi_{l,2}, ..., \xi_{l,N})$. For any measurable function $f$ such that $E|f(\bar{\xi}_1, \bar{\xi}_2, ..., \bar{\xi}_s)| < \infty$ one has the following Bartlett's type formula (see Bartlett (1938)):

$$Ef(\bar{\eta}_1, \bar{\eta}_2, ..., \bar{\eta}_s)$$

$$= \frac{1}{\prod_{l=1}^{s}(2\pi P\{\zeta_{l,N} = n_l\})} \int_{-\pi}^{\pi} ... \int_{-\pi}^{\pi} Ef(\bar{\xi}_1, \bar{\xi}_2, ..., \bar{\xi}_s) \exp\left\{i\sum_{l=1}^{s}\tau_l(\zeta_{l,N} - n_l)\right\} d\tau_1...d\tau_s. \tag{2.3}$$

Indeed, by the total expectation formula

$$Ef(\bar{\xi}_1, \bar{\xi}_2, ..., \bar{\xi}_s)\exp\left\{i\sum_{l=1}^{s}\tau_l(\zeta_{l,N} - n_l)\right\} =$$

$$E_{\zeta_{1,N}, \zeta_{2,N}, ..., \zeta_{s,N}}\left[\exp\left\{i\sum_{l=1}^{s}\tau_l(\zeta_{l,N} - n_l)\right\} E\left(f(\bar{\xi}_1, \bar{\xi}_2, ..., \bar{\xi}_s)/\zeta_{1,N}, \zeta_{2,N}, ..., \zeta_{s,N}\right)\right]$$

Integrating both sides of this equality with respect to $\tau_1, \tau_2, ..., \tau_s$ over the $s$ dimensional cube $[-\pi, \pi]^s$ we get

$$\int_{-\pi}^{\pi} \ldots \int_{-\pi}^{\pi} Ef(\bar{\xi}_1, \bar{\xi}_2, \ldots, \bar{\xi}_s) \exp\left\{i \sum_{l=1}^{s} \tau_l (\zeta_{l,N} - n_l)\right\} d\tau_1 d\tau_2 \ldots d\tau_s$$

$$= (2\pi)^s P\{\zeta_{1,N} = n_1, \zeta_{2,N} = n_2, \ldots, \zeta_{s,N} = n_s\} E\left(f(\bar{\xi}_1, \bar{\xi}_2, \ldots, \bar{\xi}_s) / \zeta_{1,N} = n_1, \zeta_{2,N} = n_2, \ldots, \zeta_{s,N} = n_s\right).$$

Formula (2.3) follows due to (2.2) and the fact that $\zeta_{1,N}, \zeta_{2,N}, \ldots, \zeta_{s,N}$ are independent r.v.'s.

Bartlett's type formula seems to be very effective in applications to problems of random allocation schemes see, for instance, Holst (1979) and Mirakhmedov (1985, 1996, 2007). In particular, if $f(\bar{\eta}_1, \bar{\eta}_2, \ldots, \bar{\eta}_s) = \exp\{it\mu_0\}$ formula (2.3) generates an integral representation of the characteristic function of $\mu_0$. It is important that in this case the integrand in (2.3) is the characteristic function of a sum of independent $(s+1)$-dimensional random vectors (r.vec's) (see relation (3.10) below). This fact allows us to use the method of characteristic functions which is well developed for a sum of independent r.vec's, see e.g. the book by Bhattacharya and Rao (1976) (from now on it is referred to as BR).

## 3. Asymptotic expansion of the characteristic function of $\mu_0$.

Put $p_l = n_l / N$, $q_l = 1 - p_l$, $Q_s = q_1 q_2 \ldots q_s$,

$$g(u_1, u_2, \ldots, u_s) = \mathbb{1}\{u_1 + u_2 + \ldots + u_s = 0\} - Q_s + Q_s \sum_{l=1}^{s} q_l^{-1}(u_l - p_l), \quad (3.1)$$

and define an r.v. $\nu_0$ by

$$\nu_0 = \sum_{m=1}^{N} g(\eta_{1m}, \eta_{2m}, \ldots, \eta_{sm}). \quad (3.2)$$

We have $\nu_0 = \mu_0 - NQ_s$ and

$$E\mu_0 = \sum_{m=1}^{N} P\{\eta_{1m} + \eta_{2m} + \ldots + \eta_{sm} = 0\} = \sum_{m=1}^{N} \prod_{l=1}^{s} P\{\eta_{lm} = 0\} = NQ_s. \quad (3.3)$$

Note that $g(\xi_{1m}, \xi_{2m}, \ldots, \xi_{sm})$ as well as $(g(\xi_{1m}, \xi_{2m}, \ldots, \xi_{sm}), \xi_{1m}, \xi_{2m}, \ldots, \xi_{sm})$, $m = 1, 2, \ldots, N$, are i.i.d. r.vec's. Moreover, it is not hard to check that

$$\sum_{m=1}^{N} Var\, g(\xi_{1m}, \xi_{2m}, \ldots, \xi_{sm}) = N\sigma^2 \quad (3.4)$$

with

$$\sigma^2 = Q_s \left(1 - Q_s (1 + \sum_{l=1}^{s} p_l q_l^{-1})\right), \quad (3.5)$$

and also that

$$Eg(\xi_{1m}, \xi_{2m}, \ldots, \xi_{sm}) = 0, \quad Eg(\xi_{1m}, \xi_{2m}, \ldots, \xi_{sm})(\xi_{lm} - p_l) = 0, \quad l = 1, 2, \ldots, s, \quad (3.6)$$

$$Eg(\xi_{1m},\xi_{2m},...,\xi_{sm})(\xi_{lm}-p_l)^2=0, \ l=1,2,...,s. \tag{3.7}$$

For simplicity of notation we put

$$\tilde{g}_m = g(\xi_{1m},\xi_{2m},...,\xi_{sm})/\sigma, \qquad \tilde{\xi}_{lm} = (\xi_{lm}-p_l)/\sqrt{p_l q_l}. \tag{3.8}$$

Let $\tilde{g}$ and $\tilde{\xi}_l$ be r.v.'s having common distributions with $\tilde{g}_m$ and $\tilde{\xi}_{lm}$ respectively. Set

$$\Theta_N(t) = \int_{A_0} \Psi^N(t,\tau_1,...,\tau_s) d\tau_1...d\tau_1 \tag{3.9}$$

with

$$A_0 = \left\{\tau_1,\tau_2,...,\tau_s : |\tau_l| \le \pi\sqrt{Np_l q_l}, l=1,2,...,s\right\},$$

$$\Psi(t,\tau_1,...,\tau_s) = E\exp\left\{\frac{it}{\sqrt{N}}\tilde{g} + \sum_{l=1}^{s}\frac{i\tau_l}{\sqrt{N}}\tilde{\xi}_l\right\}.$$

Put $f(\overline{\eta}_1,\overline{\eta}_2,...,\overline{\eta}_s) = \exp\{itv_0/\sigma\sqrt{N}\}$ in the formula (2.3) and use the inversion formula for the local probability $P\{\zeta_{lN}=n_l\}$ to get

$$\varphi_N(t) \stackrel{def}{=} E\exp\left\{it\frac{v_0}{\sigma\sqrt{N}}\right\} = \frac{\Theta_N(t)}{\Theta_N(0)}. \tag{3.10}$$

Define polynomials $G_1(t)$ and $G_2(t)$ as

$$G_k(t) = \frac{1}{(2\pi)^{s/2}}\int_{R^s} P_k(t,\tau_1,...,\tau_s)\exp\left\{-\frac{1}{2}\sum_{l=1}^{s}\tau_l^2\right\}d\tau_1...d\tau_s, \ k=1,2 \tag{3.11}$$

where

$$P_1(t,\tau_1,...,\tau_s) = \frac{i^3}{6}E\left(t\tilde{g}_m + \tau_1\tilde{\xi}_{1m} + ... + \tau_s\tilde{\xi}_{sm}\right)^3,$$

$$P_2(t,\tau_1,...,\tau_s) = \frac{i^4}{24}\left[E\left(t\tilde{g}_m + \tau_1\tilde{\xi}_{1m} + ... + \tau_s\tilde{\xi}_{sm}\right)^4 - 3\left(E\left(t\tilde{g}_m + \tau_1\tilde{\xi}_{1m} + ... + \tau_s\tilde{\xi}_{sm}\right)^2\right)^2\right]$$

$$+\frac{1}{2}P_1^2(t,\tau_1,...,\tau_s).$$

Taking into account (3.6) and (3.7) we obtain

$$G_1(t) = \frac{(it)^3}{6}E\tilde{g}^3,$$

$$G_2(t) = \frac{(it)^6}{72}(E\tilde{g}^3)^2 + \frac{(it)^4}{24}\left[E\tilde{g}^4 - 3\sum_{l=1}^{s}(E\tilde{g}^2\tilde{\xi}_l)^2 - 3\right] + \frac{(it)^2}{4}\sum_{l=1}^{s}\left(E\tilde{g}^2\tilde{\xi}_l E\tilde{\xi}_l^3 - E\tilde{g}^2\tilde{\xi}_l^2 + 1\right)$$

$$+\frac{1}{24}\sum_{l=1}^{s}\left(3E\tilde{\xi}_l^4 + 5(E\tilde{\xi}_l^3)^2 + 3(s+2)\right).$$

Put

$$W_N(t) = \exp\left\{-\frac{t^2}{2}\right\}\left(1 + \frac{1}{\sqrt{N}}G_1(t) + \frac{1}{N}(G_2(t) - G_2(0))\right),$$

$$T_N = \sqrt{N}\min\left(\left(E|\tilde{g}|^3\right)^{-1}, \sqrt{p_1q_1},...,\sqrt{p_sq_s}\right). \tag{3.12}$$

**Theorem 1**. There exist constants $C_0$ and $C_1$ such that if

$$|t| \leq C_0 T_N \tag{3.13}$$

then

$$|\varphi_N(t) - W_N(t)| \leq C_1 L_N (1 + |t|^5)\exp\left\{-\frac{t^2}{12}\right\}$$

with

$$L_N = N^{-3/2}\left(E|\tilde{g}|^5 + \sum_{l=1}^{s}(p_l q_l)^{-3/2}(p_l^4 + q_l^4)\right). \tag{3.14}$$

**Proof.** Let

$$\mathbb{P}_\mathbb{N}(t,\tau_1,...,\tau_s) = \exp\left\{-\frac{1}{2}\left(t^2 + \sum_{l=1}^{s}\tau_l^2\right)\right\}\left(1 + \frac{1}{\sqrt{N}}P_1(t,\tau_1,...,\tau_s) + \frac{1}{N}P_2(t,\tau_1,...,\tau_s)\right),$$

$$A_1 = \left\{\tau_1,...,\tau_s : |\tau_l| \leq C_2\sqrt{Np_l q_l}\left(p_l^4 + q_l^4\right)^{-1/3}, l = 1, 2,..., s\right\}.$$

**Lemma 1.** If $|t| \leq \sqrt{2N}\left(E|\tilde{g}|^3\right)^{-1}$ and $(\tau_1,...,\tau_s) \in A_1$ then

$$|\Psi^N(t,\tau_1,...,\tau_s)| \leq \exp\left\{-\frac{1}{6}(t^2 + \tau_1^2 + ... + \tau_s^2)\right\}.$$

**Proof.** In view of (3.6), the r.vec. $V_m = \left(\tilde{g}_m, \tilde{\xi}_{1m},...,\tilde{\xi}_{sm}\right)$ has a unit correlation matrix. Due to the well-known inequalities between Lyapunov's ratios (see Lemma 6.2 of BR) we have

$E|\tilde{\xi}_l|^3 \leq \left(E|\tilde{\xi}_l|^5\right)^{1/3}$ and $E|\tilde{\xi}_l|^5 = (p_l^4 + q_l^4)(p_l q_l)^{-3/2}$. Now Lemma 1 follows from Theorem 8.7 of BR.

**Lemma 2**. There exist constants $C_2$ and $C_3$ such that if

$$|t| \leq C_2\sqrt{N}\left(E|\tilde{g}|^5\right)^{-1/3} \tag{3.15}$$

and $(\tau_1,...,\tau_s) \in A_1$ then

$$|\Psi^N(t,\tau_1,...,\tau_s) - \mathbb{P}_\mathbb{N}(t,\tau_1,...,\tau_s)| \leq C_3 L_N\left(1 + |t|^9 + \sum_{l=1}^{s}|\tau_l|^9\right)\exp\left\{-\frac{1}{4}\left(t^2 + \sum_{l=1}^{s}\tau_l^2\right)\right\}.$$

**Proof.** Lemma 2 follows from Theorem 9.10 of BR because (3.6) and $V_m$ has unit correlation matrix.

Recalling (3.9) and (3.11) we have

$$\Theta_N(t) - \left(1 + \frac{1}{\sqrt{N}} G_1(t) + \frac{1}{N} G_2(t)\right) \exp\left\{-\frac{t^2}{2}\right\} = \int_{A_1} \left(\Psi^N(t, \tau_1, \ldots, \tau_s) - \mathbb{P}_\mathbb{N}(t, \tau_1, \ldots, \tau_s)\right) d\tau_1 \ldots d\tau_s$$

$$- \int_{\mathbb{R}^s - A_1} \mathbb{P}_\mathbb{N}(t, \tau_1, \ldots, \tau_s) d\tau_1 \ldots d\tau_s + \int_{A_0 - A_1} \Psi^N(t, \tau_1, \ldots, \tau_s) d\tau_1 \ldots d\tau_s \stackrel{def}{=} \nabla_1 + \nabla_2 + \nabla_3. \qquad (3.16)$$

Let $t$ satisfy (3.15). Then, using Lemma 2, we have

$$|\nabla_1| \leq C_4 L_N (1 + |t|^9) \exp\left\{-\frac{t^2}{4}\right\}. \qquad (3.17)$$

Now assume that $C_2 \sqrt{N} \left(E|\tilde{g}|^5\right)^{-1/3} \leq |t| \leq C_0 \sqrt{N} \left(E|\tilde{g}|^3\right)^{-1}$. Applying Lemma 1 we obtain

$$|\nabla_1| \leq \int_{A_1} |\Psi^N(t, \tau_1, \ldots, \tau_s)| d\tau_1 \ldots d\tau_s + \int_{A_1} |\mathbb{P}_\mathbb{N}(t, \tau_1, \ldots, \tau_s)| d\tau_1 \ldots d\tau_s \leq C_5 L_N \exp\left\{-\frac{t^2}{12}\right\}. \qquad (3.18)$$

Due to the definition of $\mathbb{P}_\mathbb{N}(t, \tau_1, \ldots, \tau_s)$, it is clear that

$$|\nabla_2| \leq C_6 L_N (1 + t^6) \exp\left\{-\frac{t^2}{2}\right\}. \qquad (3.19)$$

We have

$$|\Psi(t, \tau_1, \ldots, \tau_s)| \leq \left| E \exp\left\{\frac{i}{\sqrt{N}} \sum_{l=1}^s \tau_l \tilde{\xi}_l\right\} + E\left(\left(\exp\left\{\frac{it}{\sqrt{N}} \tilde{g}\right\} - 1\right) \exp\left\{\frac{i}{\sqrt{N}} \sum_{l=1}^s \tau_l \tilde{\xi}_l\right\}\right)\right|$$

$$\leq \left| E \exp\left\{\frac{i}{\sqrt{N}} \sum_{l=1}^s \tau_l \tilde{\xi}_l\right\}\right| + \frac{|t|}{\sqrt{N}} E|\tilde{g}| = \prod_{l=1}^s \left| E \exp\left\{\frac{i\tau_l}{\sqrt{N}} \tilde{\xi}_l\right\}\right| + \frac{|t|}{\sqrt{N}} E|\tilde{g}|.$$

Use here the inequalities $x \leq \exp\{(x^2 - 1)/2\}$ and $1 + x \leq e^x$ to get

$$|\Psi(t, \tau_1, \ldots, \tau_s)| \leq \exp\left\{-\frac{1}{2} \sum_{l=1}^s \left(1 - \left|E \exp\left\{\frac{i\tau_l}{\sqrt{N}} \tilde{\xi}_l\right\}\right|^2\right)\right\} + \frac{|t|}{\sqrt{N}} E|\tilde{g}|$$

$$\leq \exp\left\{-\frac{1}{2} \sum_{l=1}^s \left(1 - \left|E \exp\left\{\frac{i\tau_l}{\sqrt{N}} \tilde{\xi}_l\right\}\right|^2\right) + \frac{e^s |t|}{\sqrt{N}} E|\tilde{g}|\right\}. \qquad (3.20)$$

Because $\tilde{\xi}_l \sim B(1, p_l)$ and

$$\sin^2 \frac{z}{2} \geq \frac{z^2}{\pi^2}, \quad |z| \leq \pi, \qquad (3.21)$$

we have

$$\left| E\exp\left\{ \frac{i\tau_l}{\sqrt{N}}\tilde{\xi}_l \right\} \right|^2 = 1 - p_l q_l \sin^2\frac{\tau_l}{2\sqrt{Np_l q_l}} \leq 1 - \frac{\tau_l^2}{\pi^2 N}.$$

Hence from (3.20) we obtain

$$\left| \Psi(t,\tau_1,...,\tau_s) \right| \leq \exp\left\{ -\frac{1}{2\pi^2 N}\sum_{l=1}^{s}\tau_l^2 + \frac{e^s |t| E|\tilde{g}|}{\sqrt{N}} \right\}. \tag{3.22}$$

If $(\tau_1,...,\tau_s) \in A_0 - A_1$ then there exists an index $l_0$ such that $C_2\sqrt{Np_{l_0}q_{l_0}}\left(p_{l_0}^4 + q_{l_0}^4\right)^{-1/3} \leq |\tau_{l_0}| \leq \pi\sqrt{Np_{l_0}q_{l_0}}$. Therefore, using (3.22), we get

$$|\nabla_3| \leq 2^{(s+1)/2}\pi^s \exp\left\{ -Np_{l_0}q_{l_0}\left[\frac{C_2}{4\pi^2} - C_0 e^s\right] \right\} \leq 2^{(s+1)/2}\pi^s \exp\left\{ -\frac{C_2}{8\pi^2}Np_{l_0}q_{l_0} \right\} \tag{3.23}$$

because of (3.13), (3.12), the choice of $C_0 \leq C_2/e^s 8\pi^2$ and $E|\tilde{g}| \leq 1$ due to the well-known inequalities for moments. Thus by (3.17), (3.18), (3.19) and (3.23) from (3.16) we have

$$\left| \Theta_N(t) - \left(1 + \frac{1}{\sqrt{N}}G_1(t) + \frac{1}{N}G_2(t)\right)\exp\left\{-\frac{t^2}{2}\right\} \right| \leq C_7 L_N \left(1+|t|^9\right)\exp\left\{-\frac{t^2}{12}\right\}.$$

In particular,

$$\left| \Theta_N(0) - \left(1 + \frac{1}{N}G_2(0)\right) \right| \leq C_7 L_N. \tag{3.24}$$

Theorem 1 follows from these relations and formula (3.10).

**Remark 1**. Let, for some $C$,

$$E|\tilde{g}|\left(E|\tilde{g}|^3\right)^{-1} \leq C\min_{1\leq l\leq s}(p_l g_l). \tag{3.25}$$

Then condition (3.13) of Theorem 1 can be replaced by

$$|t| \leq C_0 \sqrt{N}\left(E|\tilde{g}|^3\right)^{-1}. \tag{3.26}$$

Indeed, in this case we should only give an alternative bound for $\nabla_3$ as follows:

$$|\nabla_3| \leq 2^{(s+1)/2}\pi^s \exp\left\{ -Np_{l_0}q_{l_0}\left[\frac{C_2}{4\pi^2} - C_0 e^s \frac{E|\tilde{g}|}{p_{l_0}q_{l_0}E|\tilde{g}|^3}\right] \right\}$$

$$\leq 2^{(s+1)/2}\pi^s \exp\left\{ -Np_{l_0}q_{l_0}\left[\frac{C_2}{4\pi^2} - CC_0 e^s\right] \right\} \leq 2^{(s+1)/2}\pi^s \exp\left\{ -\frac{C_2}{8\pi^2}Np_{l_0}q_{l_0} \right\},$$

because of (3.25), (3.26) and due to the choice of $C_0 \leq C_2/Ce^s 8\pi^2$.

## 4. Asymptotic Expansion in the Local Limit Theorem for $\mu_0$. Case of fixed s.

We set

$$M_2 = \sum_{l=1}^{s}\left(E\tilde{g}^2\tilde{\xi}_l E\tilde{\xi}_l^3 - E\tilde{g}^2\tilde{\xi}_l^2 + 1\right), \quad M_3 = E\tilde{g}^3, \quad M_4 = E\tilde{g}^4 - 3\sum_{l=1}^{s}\left(\left(E\tilde{g}^2\tilde{\xi}_l\right)^2 + 1\right)$$

and define $\hat{W}_N(x)$ as the inverse Fourier transform of above defined function $W_N(t)$. The function $\hat{W}_N(x)$ can be obtained by formally substituting $(-1)^\nu \dfrac{d^\nu}{dx^\nu}e^{-x^2/2}$ instead of $(it)^\nu e^{-t^2/2}$ for each $\nu$ in the expression for $W_N(t)$. Let $H_\nu(x)$ stand for the $\nu$-th order Hermit-Chebishev polynomial. Recall that

$$H_2(x) = x^2 - 1, \; H_3(x) = x^3 - 3x, \; H_4(x) = x^4 - 6x^2 + 3, \; H_6(x) = x^6 - 15x^4 + 45x^2 - 15.$$

We have

$$\hat{W}_N(x) = \frac{1}{\sqrt{2\pi}}e^{-\frac{x^2}{2}}\left(1 + \frac{H_3(x)}{6\sqrt{N}}M_3 + \frac{1}{4N}\left(\frac{1}{18}H_6(x)M_3^2 + \frac{1}{8}H_4(x)M_4 + H_2(x)M_2\right)\right).$$

Recall (3.1), (3.5), (3.8) and additionally put

$$n = \max_{1\le l\le s} n_l, \quad P_s = p_1 \cdot \ldots \cdot p_s, \quad \alpha = \max_{1\le l\le s}\min\left(Q_s^2 p_l q_l^{-1}, P_s^2 p_l^{-1} q_l\right),$$

$$x_k = (k - NQ_s)/\sqrt{N}\sigma. \tag{4.1}$$

Also keep in mind that all of $p_l, \sigma, \tilde{g}$ etc depend on $N, n_1, \ldots, n_s$.

**Theorem 2**. Let (3.25) hold and

$$\sup_{n_1,\ldots,n_s}\sigma E|\tilde{g}|^3 \le C. \tag{4.2}$$

Then there exist constants $C_8$ and $C_9$ such that

$$\max_{0\le k\le N-n}\left|\sqrt{N}\sigma P\{\mu_0 = k\} - \hat{W}_N(x_k)\right| \le C_8\left(L_N + N^{(s+1)/2}\sigma P_s Q_s e^{-C_9 N\alpha}\right),$$

with $L_N$ and $\sigma$ from (3.12) and (3.5) respectively.

From Theorem 2 immediately follows

**Corollary 1**. If $p_l$ are bounded away from zero and one for all $l = 1,2,\ldots,s$ then

$$\max_{0\le k\le N-n}\left|\sqrt{N}\sigma P\{\mu_0 = k\} - \hat{W}_N(x_k)\right| = O\left(\frac{1}{N^{3/2}}\right).$$

**Proof of Theorem 2**. Put $T_N = C_0\sqrt{N}\left(E|\tilde{g}|^3\right)^{-1}$. Due to the inversion formula, we have

$$\left|\sqrt{N}\sigma P\{\mu_0 = k\} - \hat{W}_N(x_k)\right| \le \frac{1}{\sqrt{2\pi}}\int_{|t|\le T_N}|\varphi_N(t) - W_N(t)|\,dt + \int_{T_N\le|t|\le\pi\sqrt{N}\sigma}|\varphi_N(t)|\,dt$$

$$+ \int_{|t| \geq T_N} |W_N(t)| dt \stackrel{def}{=} \Delta_1 + \Delta_2 + \Delta_3. \tag{4.3}$$

Use Theorem 1 and Remark 1 to get

$$\Delta_1 \leq C_{10} L_N. \tag{4.4}$$

It is obvious that

$$\Delta_3 \leq C_{11} L_N. \tag{4.5}$$

Let $X^*$ be the symmetrization of the r.vec. X, i.e. $X^* = X - X'$, where $X'$ and $X$ are mutually independent r.vec.'s having a common distribution, $\langle a \rangle$ be the distance between the real $a$ and the integers; $D_X(u) = E \langle (X^*, u) \rangle^2$, where $(X^*, u)$ is the scalar product of the vectors $X^*$ and $u$.

**Lemma 3.** For any real $u$ one has

$$4 D_X \left( \frac{u}{2\pi} \right) \leq 1 - |E \exp\{i(X, u)\}|^2 \leq 2\pi^2 D_X \left( \frac{u}{2\pi} \right).$$

**Proof** of Lemma 3 was presented by Mukhin (1984).

For simplicity of notation, let $g$ and $\xi_l$ stand for r.v.'s having the same distributions as $g(\xi_{1m}, \ldots, \xi_{sm})$ and $\xi_{lm}$, respectively.

Let $C_{12} \leq |t| \leq \pi$. Putting $X = (g, \xi_1, \ldots, \xi_s)$ and $u = (t, \tau_1, \ldots, \tau_s)$, we have

$$D_X \left( \frac{u}{2\pi} \right) = E \left\langle \frac{1}{2\pi} \left( g^* t + \xi_1^* \tau_1 + \ldots + \xi_s^* \tau_s \right) \right\rangle^2$$

$$\geq P\{\xi_1 + \ldots + \xi_s = 0, \xi_1' = 1, \xi_2' + \ldots + \xi_s' = 0\} \left\langle \frac{1}{2\pi} \left[ (1 + Q_s q_1^{-1}) t - \tau_1 \right] \right\rangle^2$$

$$+ P\{\xi_1 = 0, \xi_2 = \ldots = \xi_s = \xi_1' = \xi_2' = \ldots = \xi_s' = 1\} \left\langle \frac{1}{2\pi} \left[ Q_s q_1^{-1} t - \tau_1 \right] \right\rangle^2$$

$$= Q_s^2 p_1 q_1^{-1} \left\langle \frac{1}{2\pi} \left[ (1 + Q_s q_1^{-1}) t - \tau_1 \right] \right\rangle^2 + P_s^2 q_1 p_1^{-1} \left\langle \frac{1}{2\pi} \left[ Q_s q_1^{-1} t - \tau_1 \right] \right\rangle^2 \geq C_{13} \alpha, \tag{4.6}$$

because $t$ is bounded away from zero and hence $\left\langle \frac{1}{2\pi} \left[ (1 + Q_s q_1^{-1}) t - \tau_1 \right] \right\rangle^2$ and $\left\langle \frac{1}{2\pi} \left[ Q_s q_1^{-1} t - \tau_1 \right] \right\rangle^2$ can not be equal to zero simultaneously. Now use inequalities $x \leq \exp\{(x^2 - 1)/2\}$, (4.6) and Lemma 3 to get

$$\left| \Psi \left( t\sigma \sqrt{N}, \tau_1 \sqrt{N p_1 q_1}, \ldots, \tau_s \sqrt{N p_s q_s} \right) \right|$$

$$\leq \exp\left\{-\frac{1}{2}\left(1-\left|\Psi\left(t\sigma\sqrt{N},\tau_1\sqrt{Np_1q_1},...,\tau_s\sqrt{Np_sq_s}\right)\right|^2\right)\right\} \leq e^{-C_{14}\alpha/2}. \qquad (4.7)$$

Hence using (3.10) and taking into account (3.9), (3.10), (3.24) and condition (4.2), under which $T_N/\sigma\sqrt{N} \geq C_0/C = C_{12}$, by simple manipulations we obtain

$$|\nabla_2| \leq \frac{\pi^{s+1}}{\Theta(0)}\sqrt{P_sQ_s}N^{(s+1)/2}\exp\{-C_{14}N\alpha/2\} \leq \exp\{-C_{15}N\alpha\}. \qquad (4.8)$$

Theorem 2 follows from (4.3),(4.4), (4.5) and (4.8).

## 5. Asymptotic Expansion in the Local Limit Theorem for $\mu_0$ when s may increase.

In this Section we apply another approach (known as Frobenius-Harper technique) first considered by Frobenius (1910) and rediscovered by Harper (1967); it was also used, for instance, by Park (1981), Vatutin and Mikhaylov (1984) and Gani (2004). We wish to use the distributional equality of $\mu_0$ to a sum of independent Bernoulli r.v.'s. Vatutin and Mikhaylov (1984) showed that the probability generating function $F(z) = Ez^{\mu_0}$ satisfies the Frobenius-Harper property: $F(z)$ is a polynomial of degree $N-n$ and has only negative real roots, which we denote by $-d_1,-d_2,...,-d_{N-n}$. Hence

$$F(z) = \prod_{m=1}^{N-n}\frac{z-d_m}{1+d_m} = \prod_{m=1}^{N-n}Ez^{Y_m}. \qquad (5.1)$$

where $Y_1,Y_2,...,Y_{N-n}$ is a sequence of independent r.v.'s and $Y_l \sim B(1,a_l)$ with $a_l = (1+d_l)^{-1}$, $l=1,2,...,N-n$. Thus the r.v.'s $\mu_0$ and $\mu_0' = Y_1+Y_2+...+Y_{N-n}$ have a common distribution. Using this fact we can obtain an asymptotic expansion result without any restrictions on $s$ – the number of the sets of particles. Indeed, by Fourier inversion formula

$$P\{\mu_0 = k\} = \frac{1}{2\pi\sqrt{Var\mu_0}}\int_{|t|\leq\pi\sqrt{Var\mu_0}}e^{-itu_k}\gamma(t)dt,$$

where

$$u_k = (k-E\mu_0)/\sqrt{Var\mu_0}, \qquad (5.2)$$

$$\gamma(t) = E\exp\left\{it\frac{\mu_0-E\mu_0}{\sqrt{Var\mu_0}}\right\} = \prod_{m=1}^{N-n}E\exp\left\{it\frac{Y_m-a_m}{\sqrt{Var\mu_0'}}\right\}, \qquad (5.3)$$

Because of (5.1) and $E\mu_0 = E\mu_0' = a_1+a_2+...+a_{N-n}$, $Var\mu_0 = Var\mu_0'$. Due to the second equality in (5.3), for $|t| \leq \pi\sqrt{Var\mu_0'}$ we find, using (3.21), that

$$|\gamma(t)|^2 = \prod_{m=1}^{N-n}\left(1 - a_m(1-a_m)\sin^2\frac{t}{2\sqrt{Var\mu_0'}}\right) \le \exp\left\{-\frac{t^2}{\pi^2}\right\}.$$

Therefore, applying the classical method of the proof of asymptotical expansion results in the central local limit theorem for sums of independent lattice r.v.'s (see, for instance, Petrov (1995), and also Section 4 of the present paper), we can easily derive an asymptotic expansion result for $P\{\mu_0 = k\}$ without any restrictions on $s$. This result has the following form. Put

$$L_3(n,N) = (Var\mu_0')^{-3/2}\sum_{m=1}^{N-n}E(Y_m - a_m)^3 = \left(\sum_{m=1}^{N-n}a_m(1-a_m)\right)^{-3/2}\sum_{m=1}^{N-n}a_m(1-a_m)(1-2a_m),$$

$$L_4(n,N) = (Var\mu_0')^{-2}\sum_{m=1}^{N-n}\left(E(Y_m - a_m)^4 - 3\left(E(Y_m - a_m)^2\right)^2\right)$$

$$= \left(\sum_{m=1}^{N-n}a_m(1-a_m)\right)^{-2}\sum_{m=1}^{N-n}a_m(1-a_m)\left(1 - 6a_m(1-a_m)\right).$$

**Theorem 3.** There exists a constant $C$ such that

$$\max_{0 \le k \le N-n}\left|\sqrt{Var\mu_0}P\{\mu_0 = k\} - \frac{1}{\sqrt{2\pi}}e^{-\frac{u_k^2}{2}}\left(1 + \frac{H_3(u_k)}{6}L_3(n,N)\right.\right.$$

$$\left.\left. + \frac{H_6(u_k)}{72}L_3^2(n,N) + \frac{H_4(u_k)}{24}L_4(n,N)\right)\right| \le \frac{C}{(Var\mu_0)^{3/2}}.$$

However, such a result is not practical, because the terms of the asymptotic expansion depend on $d_1, d_2, ..., d_{N-n}$ the calculation of which is a difficult problem.

Theorem 2 of Section 4 assumed that the number $s$ of the sets of particles is fixed. This is a technical condition, as we used the method of asymptotic expansions for the characteristic functions of the sums of independent r. vec.'s, where $s$ is the dimension of those vectors. Actually, the basic formula (3.10) is still correct without any restrictions on $s$. Hence it can be used for formal construction of the terms of an asymptotic expansion of the characteristic function of $(\mu_0 - E\mu_0)/\sqrt{Var\mu_0}$.

Recalling (2.1) and that $\eta_{1m}, \eta_{2m}, ..., \eta_{sm}$ are mutually independent r.v.'s with $\eta_{lm} \sim Bi(1, p_l)$ for each $m = 1, 2, ..., N$, we have

$$E\mu_0^2 = \sum_{m=1}^{N}E\mathbb{1}\{\eta_{1m} + \eta_{2m} + ... + \eta_{sm} = 0\}$$

$$+ \sum_{\substack{m,l=1 \\ m \ne l}}^{N}E\mathbb{1}\{\eta_{1m} + \eta_{2m} + ... + \eta_{sm} = 0\}\mathbb{1}\{\eta_{1l} + \eta_{2l} + ... + \eta_{sl} = 0\}$$

$$= NQ_s + + \sum_{\substack{m,l=1 \\ m \neq l}}^{N} P\{\eta_{1m} = 0, \eta_{1l} = 0, \eta_{2m} = 0, \eta_{2l} = 0, ..., \eta_{sm} = 0, \eta_{sl} = 0\}.$$

The r.vec.'s $(\eta_{1m}, \eta_{1l}), (\eta_{2m}, \eta_{2l}), ..., (\eta_{sm}, \eta_{sl})$ are mutually independent, because the sets of particles are allocated independent of each other. Hence

$$E\mu_0^2 = NQ_s + \sum_{\substack{m,l=1 \\ m \neq l}}^{N} \prod_{j=1}^{s} P\{\eta_{jm} = 0, \eta_{jl} = 0\} = NQ_s + N(N-1)\prod_{j=1}^{s}\left(1 - \frac{n_j}{N}\right)\left(1 - \frac{n_j}{N-1}\right)$$

$$= NQ_s + N(N-1)Q_s^2 - NQ_s^2 \sum_{i=1}^{s} p_i q_i^{-1} + NQ_s^2 \sum_{\nu=2}^{s} \frac{(-1)^\nu}{(N-1)^{\nu-1}} \sum p_{i_1} p_{i_2} \cdots p_{i_\nu} (q_{i_1} q_{i_2} \cdots q_{i_\nu})^{-1}.$$

where $\sum$ denotes summation over all $\nu$-tuples $(i_1, ..., i_\nu)$ with components not equal to each other and such that $i_m = 1, 2, ..., s$ ; $m = 1, 2, ..., \nu$. From this and (3.3), (3.5) we have

$$Var\mu_0 = N\sigma^2 + \frac{N}{N-1}Q_s^2 \sum_{\nu=2}^{s} \frac{(-1)^{\nu-2}}{(N-1)^{\nu-2}} \sum\nolimits^* p_{i_1} p_{i_2} \cdots p_{i_\nu} (q_{i_1} q_{i_2} \cdots q_{i_\nu})^{-1}$$

$$= N\sigma^2 \left(1 + \frac{Q_s}{(N-1)\left(1 - Q_s\left(1 + \sum_{l=1}^{s} p_l q_l^{-1}\right)\right)} \sum_{\nu=2}^{s} \frac{(-1)^{\nu-2}}{(N-1)^{\nu-2}} \sum p_{i_1} p_{i_2} \cdots p_{i_\nu} (q_{i_1} q_{i_2} \cdots q_{i_\nu})^{-1}\right)$$

$$\stackrel{def}{=} N\sigma^2(1 + b_N). \tag{5.4}$$

Note that $\gamma(t) = \varphi_N\left((N\sigma^2/Var\mu_0)^{1/2} t\right) = \varphi_N\left((1+b_N)^{-1/2} t\right)$ (see (3.10) and (5.3)) owing to (3.4) and (5.4). So due to (3.10) we formally have $\gamma(t) = W_N\left((1+b_N)^{-1/2} t\right) + \varepsilon_N(t) = W_N(t) + \frac{(it)^2}{2} b_N + \tilde{\varepsilon}_N(t)$, where $\varepsilon_N(t)$ and $\tilde{\varepsilon}_N(t)$ are of the order $O(N^{-3/2})$. Hence the first three terms of asymptotic expansion of $\sqrt{Var\mu_0} P(\mu_0 = k)$ equal to $\hat{W}_N(u_k) + N^{-1} H_2(u_k) \tilde{b}_N$, where $\tilde{b}_N = Nb_N$. This, in combination with Theorem 3, implies the following result

**Theorem 4**. There exists a constant $C$ such that

$$\max_{0 \leq k \leq N-n} \left| \sqrt{Var\mu_0} P\{\mu_0 = k\} - \hat{W}_N(u_k) - \frac{1}{N} H_2(u_k) \tilde{b}_N \right| \leq \frac{C}{\left(N\sigma^2(1+b_N)\right)^{3/2}}.$$

**Remark 2**. In Theorem 4 exact normalization by $\sqrt{Var\mu_0}$ is used, see $u_k$ in (5.2), whereas in Theorem 2 we deal with normalization by $\sigma\sqrt{N}$ ( see $x_k$ in (4.1)), which is an asymptotic of $\sqrt{Var\mu_0}$. The extra term $N^{-1}H_2(u_k)\tilde{b}_N$ appeared in Theorem 4 for that reason.

**Remark 3**. The r.v. $\mu_0$ is the special case of the statistics of the form $\sum_{m=1}^{N}h(\eta_{1,m},\eta_{2,m},...,\eta_{s,m})$, were $h$ is a measurable real function. For this general statistics Bartlett type formula is also hold. Hence the formal construction of the terms of asymptotic expansions can be established. Some applications of Bartlett type formula for so called "generalized allocation schemes" are presented in Mirakhmedov (1985, 1994, 1995, 1996, 2007) .

**Remark 4**. We restricted ourselves to the first three terms of the asymptotic expansion. It is obvious that the presented approach can be used to derive asymptotic expansions of any length.

**Acknowledgement**.
We would like to express our gratitude to the referees for useful comments. We also thankful to Prof. Kostya Borovkov for careful reading of our manuscript and help us to improve the presentation of the paper.

**References**.
1. Bartlett M.S., 1938. The characteristic function of a conditional statistics. J. London Math., 13, p. 62-67.
2. Bhattachatya R.N. and R. Ranga Rao, 1976. Normal approximation and Asymptotic Expansions. John Wiley&Sons, New-York.
3. Feller, W., 1968. An introduction to Probability Theory and Its Applications. V.1, third ed. Wiley, New York.
4. Frobenius G., 1910. Über die Bernoullischen Zahlen und die Eulerschen Polynome. Sitzenberichte der Preussischen Akademie der Wissenschaften, p. 829-830
5. Gani, J., 1993. Random allocation methods in an epidemic model. In Stochastic Processes: A Festachrift in Honour of Gopinath Kallianpur, S.Cambanis,J.K.Ghosh, R.L.Karandikar and P.K.Sen (Editors).Springer. New-York. 97-106.
6. Gani, J., 2004.  Random allocation and urn models. J.Appl.Probabl., 41A, p.313-320.
7. Ivanov V.A., Ivchenko G.I., and Medvedev Yu.I., 1985. Discrete problems of the probability theory (a survey).J. Soviet Math.,31, no.2 .
8. Harper L., 1967. Stirling behavior is asymptotically normal. Ann. Math. Statist. 38, p. 410-414.

9. Holst L., 1979. A unified approach to limit theorems for urn models. J.Appl., Probab., 16,#1,p. 154-162.
10. Jonson N.L., Kotz S., 1977. Urn Models and their Applications. Wiley, New York.
11. Kolchin V.F., Sevastyanov B.A. and Chistyakov V.P., 1978. Random Allocation. John Wiley, New-York-Toronto-London.
12. Kotz S. and Balakrishan N., 1997. Advances in urn models during the past two decades. In Advances in Combinatorial Methods and Appl. to Probabl. and Statist., p.203-257. Birkhauser, Boston. MA.
13. Mikhaylov V.G., 1978. Asymptotic normality of the number of empty cells for group allocation of particles. Theory Probabl. Appl., 25, p.82-90.
14. Mikhaylov V.G., 1981.Convergence to multi-dimensional normal law in an equiprobable scheme for group allocation of particles. Math USSR-St, 39, p.145-168.
15. Menezes, A., Van Qorshot, P. and S. Vanstone, 1997. Handbook of applied cryptology. CRC Press, New York.
16. Mirakhmedov Sh. A., 1985. The estimations of the closeness to a normal law in the scheme without replacement. Theory Probabl.Appl.,v.30, p.427-439.
17. Mirakhmedov Sh. A., 1994. Asymptotical analysis of the conditional distributions of the sum of independent random variables and applications in the theory of decomposable statistics. Dissertation for Doctor of Sciences degree. Tashkent. Uzbekistan. 290 p.
18. Mirakhmedov Sh. A., 1995, Limit theorems for conditional distributions. Discrete Math. Appl., v.5, p. 107-132.
19. Mirakhmedov Sh.A., 1996. Limit theorems on decomposable statistics in a generalized allocation schemes. Discrete Math. Appl., v. 6, p. 379-404.
20. Mirakhmedov Sh.M., 2007. Asymptotic normality associated with generalized occupancy problem. Statist. & Probabl. Letters, v.77, p. 1549-1558.
21. Mukhin A.B., 1984. Local limit theorems for the distribution of a sum of independent random vectors. Theory Probabl. Appl., 28, p. 360-366.
22. Park C.J., 1981. On the distribution of the number of unobserved elements when m samples of size n are drawn from a finite population. Commun. Statist., A10, p.371-383.
23. Petrov V.V., 1995. Limit theorems of probability theory. Oxford University Press, New York
24. Vatutin V.A. and Mikhaylov V.G., 1982. Limit theorems for the number of empty cells in an equiprobable scheme for group allocation of particles. Theory Probabl. Appl, 27, p.734-743.